\documentclass[aos,MSNbibl,seceqn,dvips]{arximspdf}

\doi{10.1214/13-AOS1190} 
\volume{42}
\issue{1}
\pubyear{2014}
\firstpage{211}
\lastpage{224}

\makeatletter
\newcommand{\rrVert}{\Vert}
\newcommand{\llVert}{\Vert}
\newtheorem{prop}{Proposition}
\newproclaim{hyp}{Assumption}
\newtheorem{Theorem}{Theorem}
\makeatother

\begin{document}
\begin{frontmatter}

\title{Optimal learning with $Q$-aggregation}
\runtitle{Optimal learning with $Q$-aggregation}

\begin{aug}
\author[A]{\fnms{Guillaume} \snm{Lecu\'e}\thanksref{t1}\ead[label=e1]{guillaume.lecue@univ-mlv.fr}}
\and
\author[B]{\fnms{Philippe} \snm{Rigollet}\corref{}\thanksref{t2}\ead[label=e2]{rigollet@princeton.edu}}
\runauthor{G. Lecu\'e and P. Rigollet}
\affiliation{CNRS, Ecole Polytechnique and Princeton University}
\address[A]{CNRS, CMAP\\
Ecole Polytechnique\\
Palaiseau, 91120\\
France\\
\printead{e1}} 
\address[B]{Department of Operations Research\\
\quad and Financial Engineering\\
Princeton University\\
Princeton, New Jersey 08544\\
USA\\
\printead{e2}}
\end{aug}
\thankstext{t1}{Supported by French Agence Nationale de la Recherche
ANR Grant ``\textsc{Prognostic}'' ANR-09-JCJC-0101-01.}
\thankstext{t2}{Supported in part by NSF Grants DMS-09-06424,
DMS-13-17308, CAREER-DMS-10-53987, a Howard B.~Wentz Jr. award and a
gift from the Bendheim Center for Finance.}

\received{\smonth{7} \syear{2013}}
\revised{\smonth{11} \syear{2013}}

%
\begin{abstract}
We consider a general supervised learning problem with strongly convex
and Lipschitz loss and study the problem of model selection
aggregation. In particular, given a finite dictionary functions
(learners) together with the prior, we generalize the results obtained
by Dai, Rigollet and Zhang [\textit{Ann. Statist.} \textbf{40} (2012) 1878--1905]
for Gaussian regression with squared
loss and fixed design to this learning setup. Specifically, we prove
that the $Q$-aggregation procedure outputs an estimator that satisfies
optimal oracle inequalities both in expectation and with high
probability. Our proof techniques somewhat depart from traditional
proofs by making most of the standard arguments on the Laplace
transform of the empirical process to be controlled.
\end{abstract}

%
\begin{keyword}[class=AMS]
\kwd[Primary ]{68Q32}
\kwd[; secondary ]{62G08}
\kwd{62G05}
\end{keyword}
\begin{keyword}
\kwd{Learning theory}
\kwd{empirical risk minimization}
\kwd{aggregation}
\kwd{empirical processes theory}
\end{keyword}

\end{frontmatter}

\section{Introduction and main results} \label{secintroduction}
Let $\mathcal{X}$ be a probability space and let $(X,Y) \in\mathcal
{X}\times\mathbb{R}$ be a random couple. Broadly speaking, the goal of
statistical learning is to predict $Y$ given $X$. To achieve this goal,
we observe a dataset $\mathcal{D}=\{(X_1,Y_1),\ldots,(X_n,Y_n)\}$
that consists of $n$ independent copies of $(X,Y)$ and use these
observations to construct a function (learner) $f\dvtx \mathcal
{X}\rightarrow\mathbb{R}$ such that $f(X)$ is close to $Y$ in a certain
sense. More precisely, the prediction quality of a (possibly data
dependent) function $\hat f$ is measured by a \emph{risk function} $R\dvtx
\mathbb{R}^\mathcal{X}\to\mathbb{R}$ associated to a \emph
{loss function}
$\ell\dvtx  \mathbb{R}^2 \to\mathbb{R}$ in the following way:
%
\[
R(\hat f)=\mathbb{E} \bigl[\ell\bigl(Y,\hat f(X)\bigr)\vert\mathcal{D}
\bigr].
\]
We focus hereafter on loss functions $\ell$ that are \emph{convex} in
their second argument. Moreover, for the sake of simplicity, throughout
this article we restrict ourselves to functions
$f$ and random variables $(X,Y)$ for which
$|Y|\leq b$ and $|f(X)| \leq b$ almost surely, for some fixed $b\geq
0$. For any real-valued measurable $f$ on $\mathcal{X}$, for which
this quantity is finite, we define $\|f\|_2=\sqrt{\mathbb
{E}[f(X)^2]}$.

We are given a finite set $\mathcal{F}=\{f_1,\ldots,f_M\}$ of
measurable functions from $\mathcal{X}$ to~$\mathbb{R}$. This set is
called a
\emph{dictionary}. The elements in $\mathcal{F}$ may have been
constructed using an independent, frozen, dataset at some previous step
or may simply be good candidates for the learning task at hand. To
focus our contribution on the aggregation problem, we restrict our
attention to the case where $\mathcal{F}$ consists of deterministic
functions and, because of the diversity of dictionaries that can be
considered, we do not want to assume anything on the dictionary except
boundedness.

The aim of model selection aggregation \cite{Nem00,Cat04,Cat07,Tsy03}
is to use the data $\mathcal{D}$ to construct a function $\hat f$
having an \emph{excess-risk} $ R(\hat f)- \min_{f\in\mathcal
{F}}R(f)$ as small as possible. Namely, we seek the smallest
deterministic \emph{residual term} $\Delta_{n}(\mathcal{F})>0$ such
that the excess risk is bounded above by $\Delta_{n}(\mathcal{F})$,
either in expectation or with high probability, or, in this instance,
in both. In the high probability case, such bounds are called \emph
{oracle inequalities}. This problem was introduced and studied in \cite
{Cat04,Nem00}. Many results have been obtained in aggregation theory
during the last decade, for instance, in \cite{Aud08}, the
suboptimality in deviation of the Gibbs aggregates is proved, in \cite
{Aud09}, several procedures related to Gibbs aggregates are proved to
be optimal (in expectation) even under moment assumptions, in \cite
{BunTsyWeg07c}, a ``universal'' aggregation method is constructed to
solve several type of aggregation problems in the Gaussian regression
model. Other construction of optimal aggregation procedures in various
setups can also be found in
\cite{JudNem00,JudRigTsy08,LecMen09,Tsy03,Tsy04b,Yan00a,Yan00}.

From a minimax standpoint, it has been proved that $\Delta_n(\mathcal
{F})=C(\log M)/n$, \mbox{$C>0$} is the smallest residual term that one can
hope for the regression problem with quadratic loss \cite{Tsy03}. An
estimator $\hat f$ achieving such a rate (up to some multiplying
constant) is called an optimal aggregate. One of the first procedures
proved to achieve this optimal rate is a progressive mixture rule of
Gibbs estimators (cf. \mbox{\cite{Cat04,Aud09,Yan00,JudRigTsy08}}). The
optimality of this procedure holds for any ``exponentially concave''
loss function (cf. Theorem~4.2 in \cite{JudRigTsy08}) but only in
expectation (cf. \cite{Aud08}).

The aim of this paper is to construct optimal aggregates (both in
expectation and deviation) under general conditions on the loss
function $\ell$ and for a random design. We also want these procedures
to have the ability to take into account some prior information on the
dictionary unlike the existing optimal aggregation procedures that have
been constructed in this setup so far (cf. \cite{Aud08,LecMen09}).

Note that the optimal residuals for model selection aggregation are of
the order $1/n$ as opposed to the standard parametric rate $1/\sqrt
{n}$. This \emph{fast} rate essentially comes from the strong
convexity of the quadratic loss. In what follows, we show that indeed,
strong convexity is sufficient to obtain fast rates. It is known that
rates of optional order $1/n$ cannot be achieved if the loss function
is only assumed to be convex. Indeed, it follows from
\cite{Lec07a}, Theorem~2, that if the loss is linear then the best
achievable residual term is at least of the order $\sqrt{(\log
|\mathcal{F}|)/n}$. Recall that a function $g$ is said to be strongly
convex on a nonempty convex set $C \subset\mathbb{R}$ if there exists a
constant $c$ such that
\[
g\bigl(\alpha a + (1-\alpha)a'\bigr) \le\alpha g(a) + (1-\alpha) g
\bigl(a'\bigr) -\frac
{c}{2}\alpha(1-\alpha)
\bigl(a-a'\bigr)^2
\]
for any $a,a' \in C, \alpha\in(0,1)$.
In this case, $c$ is called \emph{modulus of strong convexity}.
For technical reasons, we will also need to assume that the loss
function is Lipschitz. We now introduce the set of assumptions that are
sufficient for our approach.
%
%
\begin{hyp} \label{assloss}
The loss function $\ell$ is such that for any $f,g\in[-b,b]$, we have
\[
\bigl|\ell(Y,f)-\ell(Y,g)\bigr|\leq C_b |f-g|\qquad\mbox{a.s.}
\]
Moreover, almost surely, the function $\ell(Y, \cdot)$ is \emph
{strongly convex} with modulus of strong convexity $C_\ell$ on $[-b,b]$.
\end{hyp}

A central quantity that is used for the construction of aggregates is
the empirical risk defined by
%
%
\begin{equation}
\label{eqempirical-risk} R_n(f)=\frac{1}{n}\sum
_{i=1}^n\ell\bigl(Y_i,f(X_i)
\bigr)
\end{equation}
for any real-valued function $f$ defined over $\mathcal{X}$.
A natural aggregation procedure consists in taking the function in
$\mathcal{F}$ that minimizes
the empirical risk. This procedure is called empirical risk
minimization (ERM). It has
been proved that ERM is suboptimal for the aggregation problem (cf.
Proposition~2.1 in
\cite{JudRigTsy08} or Chapter~3.5 in \cite{Cat04}, Theorem~1.1 in
\cite{LecMen10}, Theorem~3 in \cite{Lec07}, Theorem~2 in \cite
{LeeBarWil98} and Theorem~2.1 in \cite{RigTsy12}). Somehow, this
procedure does not take advantage of the convexity of the loss since
the class of functions on which the empirical risk is minimized to
construct the ERM is $\mathcal{F}$, a finite set. As it turns out, the
performance of ERM relies critically on the convexity of the class of
functions on which the empirical risk is minimized \cite
{LeeBarWil98,LecMen10}. Therefore, a natural idea is to ``improve the
geometry'' of $\mathcal{F}$ by taking its convex hull $\operatorname{conv}(\mathcal{F})$ and then by minimizing the empirical risk over it.
However, this procedure is also suboptimal \cite{LecMen09,DaiRigZha12}.
The weak point of this procedure lies in the metric
complexity of the problem: taking the convex hull of $\mathcal{F}$
indeed ``improves the geometry'' of $\mathcal{F}$ but it also
increases by too much its complexity. The complexity of the convex hull
of a set can be much larger than the complexity of the set itself and
this leads to a failure of this naive convexification trick.
Nevertheless, a compromise between geometry and complexity was stricken
in \cite{Aud08} and \cite{LecMen09} where optimal aggregates have
been successfully constructed. In~\cite{Aud08}, this improvement is
achieved by minimizing the empirical risk over a carefully chosen
star-shaped subset of the convex hull of $\mathcal{F}$. In \cite
{LecMen09}, a better geometry was achieved by taking the convex hull of
an appropriate subset of $\mathcal{F}$ and then by minimizing the
empirical risk over it.

In this paper, we show that a third procedure, called $Q$-aggregation,
and that was introduced in \cite{Rig12,DaiRigZha12} for fixed design
Gaussian regression, also leads to optimal rates of aggregation. Unlike
the above two procedures that rely on finding an appropriate constraint
for ERM, $Q$-aggregation is based on a penalization of the empirical
risk but the constraint set is kept to be the convex hull of $\mathcal
{F}$. Let $\Theta$ denote the flat simplex of $\mathbb{R}^M$ defined by
\[
\Theta= \Biggl\{(\theta_1,\ldots,\theta_M)\in
\mathbb{R}^M\dvtx \theta_j\geq0, \sum
_{j=1}^M\theta_j=1 \Biggr\}
\]
and for any $\theta\in\Theta$, define the convex combination
$f_\theta=\sum_{j=1}^M \theta_j f_j$.
For any fixed~$\nu$, the $Q$-functional is defined for any~$\theta\in
\Theta$ by
%
%
\begin{equation}
\label{eqmixture-risk} Q(\theta)=(1-\nu)R_n(f_\theta)+\nu\sum
_{j=1}^M\theta_j R_n(f_j).
\end{equation}
We keep the terminology $Q$-aggregation from \cite{DaiRigZha12} in
purpose. Indeed, $Q$ stands for \emph{quadratic} and while do not
employ a quadratic loss, we exploit strong convexity in the same manner
as in~\cite{DaiRigZha12} and \cite{Rig12}. Indeed the first term in
$Q$ acts as a regularization of the linear interpolation of the
empirical risk, and is therefore a strongly convex regularization.

We consider the following aggregation procedure. Unlike the procedures
introduced in~\cite{Aud08,LecMen09}, the $Q$-aggregation procedure
allows us to put a prior weight given by a prior probability $\pi=(\pi
_1,\ldots,\pi_M)$ on each element of the dictionary $\mathcal{F}$.
This feature turns out to be crucial for applications
\cite
{AlqLou11,DalIngTsy13,DalSal11,DalTsy07,DalTsy08,DalTsy10,DalTsy12b,RigTsy11,RigTsy12}.
Let $\beta>0$ be the \emph{temperature}
parameter and $0<\nu<1$. Consider any vector of weights $\hat\theta
\in\Theta$ defined by
%
%
\begin{equation}
\label{eqaggregation-procedure} \hat\theta\in\mathop{\operatorname{argmin}}_{\theta\in\Theta
}
\Biggl[(1-\nu)R_n(f_\theta)+\nu\sum
_{j=1}^M\theta_j R_n(f_j)-
\frac{\beta}{n}\sum_{j=1}^M
\theta_j\log\pi_j \Biggr].
\end{equation}
It comes out of our analysis that $f_{\hat\theta}$ achieves an
optimal rate of aggregation if $\beta$ satisfies
%
%
\begin{equation}
\label{eqbeta-lower-bound} \beta> \max\biggl[ \frac{8C_b^2(1-\nu
)^2}{\mu}, 4\sqrt{3}bC_b(1-
\nu), \frac{C_b \nu(\nu C_b+4\mu b)}{\mu} \biggr],
\end{equation}
where $\mu=\min(\nu, 1-\nu)(C_\ell)/10 $.\vadjust{\goodbreak}

This procedure was studied in the case of fixed design in \cite
{DaiRigZha12}, where it is shown that greedy algorithms similar to the
Frank--Wolfe algorithm, can be employed to solve the optimization
problem in~(\ref{eqaggregation-procedure}). In particular, such
algorithms can yield solutions $\hat\theta$ that are very sparse:
they can have a little as two nonzero coordinates. In this case, and
when the prior $\pi$ is uniform, this two-step procedure recovers the
STAR algorithm of Audibert~\cite{Aud08}. Furthermore, unlike the STAR
algorithm, the greedy algorithm of \cite{DaiRigZha12} was shown to (i)
allow to handle any prior $\pi$ and (ii) yield better constants as
well as better numerical performance for a larger number of iterations
(see \cite{DaiRigZha12} for more details). Similar algorithms can be
employed in the present case and it follows trivially from
\cite{DaiRigZha12}, Proposition~4.1, that $n$ iterations suffice to obtain an
optimization error of the same order as the statistical error. Going
down to two iterations as in \cite{DaiRigZha12}, Theorem~4.2, requires
a more delicate analysis, similar to the one employed in~\cite
{DaiRigZha12}, but is beyond the scope of this paper.

\renewcommand{\theTheorem}{\Alph{Theorem}}
%
\begin{Theorem}
\label{THMAINUB}
Let $\mathcal{F}$ be a finite dictionary of cardinality $M$ and
$(X,Y)$ be a~random couple of $\mathcal{X}\times\mathbb{R}$ such that
$|Y|\leq b$ and $\max_{f \in\mathcal{F}}|f(X)|\leq b$ a.s. for some
$b>0$. Assume that Assumption~\ref{assloss} holds and that $\beta$
satisfies (\ref{eqbeta-lower-bound}). Then, for any $x>0$, with
probability greater than $1-\exp(- x)$
\[
R(f_{\hat\theta})\leq\min_{j=1,\ldots,M} \biggl[R(f_j)+
\frac{\beta
}{n}\log\biggl(\frac{1}{\pi_j} \biggr) \biggr]+\frac{2\beta x}{n}.
\]
Moreover,
\[
\mathbb{E} \bigl[R(f_{\hat\theta}) \bigr]\leq\min_{j=1,\ldots,M}
\biggl[R(f_j)+\frac{\beta}{n}\log\biggl(\frac{1}{\pi_j} \biggr)
\biggr].
\]
\end{Theorem}

If $\pi$ is the uniform distribution, that is $\pi_j=1/M$ for all
$j=1, \ldots, M$, then we recover in Theorem~\ref{THMAINUB} the
classical optimal rate of aggregation $(\log M)/n$ and the estimator
$\hat\theta$ is just the one minimizing the $Q$-functional defined in
(\ref{eqmixture-risk}). In particular, no temperature parameter
$\beta$ is needed for its construction. As a result, in this case, the
parameter $b$ need not be known for the construction of the
\mbox{$Q$-}aggregation procedure.

\section{\texorpdfstring{Preliminaries to the proof of Theorem~\protect\ref{THMAINUB}}
{Preliminaries to the proof of Theorem A}}\label{secpreliminaries}

An important part of our analysis is based upon concentration
properties of empirical processes. While our proofs are similar to
those employed in \cite{Rig12} and \cite{DaiRigZha12}, they contain
genuinely new arguments. In particular, this learning setting, unlike
the denoising setting considered in \cite{Rig12,DaiRigZha12} allows
us to employ various new tools such as symmetrization and contraction.
A classical tool to quantify the concentration of measure phenomenon is
given by Bernstein's inequality for bounded variables. In terms of
Laplace transform, Bernstein's inequality \cite{BouLugMas12}, Theorem~1.10, states that if $Z_1,\ldots,Z_n$ are $n$
i.i.d. real-valued\vadjust{\goodbreak} random variables such that for all $i=1,\ldots,n$,
\[
|Z_i|\leq c\qquad\mbox{a.s.}\quad \mbox{and}\quad \mathbb{E}Z_i^2
\leq v,
\]
then for any $0<\lambda<1/c$,
%
%
\begin{equation}
\label{eqBernstein} \mathbb{E}\exp\Biggl[\lambda\Biggl(\sum
_{i=1}^n \{Z_i-\mathbb{E}Z_i
\} \Biggr) \Biggr]\leq\exp\biggl(\frac{nv\lambda^2}{2(1-c\lambda)}
\biggr).
\end{equation}
Bernstein's inequality usually yields a bound of order $\sqrt{n}$ for
the deviations of a sum around its mean. As mentioned above, such
bounds are not sufficient for our purposes and we thus consider the
following concentration result.

%
\begin{prop}\label{propBernstein-shifted}
Let $Z_1,\ldots,Z_n$ be i.i.d. real-valued random variables and let
$c_0>0$. Assume that $|Z_1|\leq c$ a.s. Then, for any $0<\lambda
<(2c_0)/(1+2c_0 c)$,
\[
\mathbb{E}\exp\Biggl[n\lambda\Biggl(\frac{1}{n}\sum
_{i=1}^nZ_i-\mathbb{E}Z_i-c_0
\mathbb{E} Z_i^2 \Biggr) \Biggr]\leq1
\]
and
\[
\mathbb{E}\exp\Biggl[n\lambda\Biggl(\frac{1}{n}\sum
_{i=1}^n\mathbb{E}Z_i-
Z_i-c_0\mathbb{E} Z_i^2
\Biggr) \Biggr]\leq1.
\]
\end{prop}
\begin{pf}
It follows from Bernstein's inequality~(\ref{eqBernstein}) that for
any $0<\lambda<(2c_0)/(1+2c_0 c)$,
\begin{eqnarray*}
&& \mathbb{E}\exp\Biggl[n\lambda \Biggl(\frac{1}{n}\sum
_{i=1}^nZ_i-\mathbb{E}Z_i-c_0
\mathbb{E} Z_i^2 \Biggr) \Biggr]
\\
&&\qquad \leq\exp\biggl(\frac{n\mathbb{E}Z_1^2\lambda^2}{2(1-c\lambda
)} \biggr) \exp\bigl[-n\lambda
c_0 \mathbb{E}Z_1^2 \bigr]
\leq1.
\end{eqnarray*}
The second inequality is obtained by replacing $Z_i$ by $-Z_i$.
\end{pf}

We also recall the following exponential bound for Rademacher
processes: let $\varepsilon_1,\ldots,\varepsilon_n$ be independent
Rademacher random variables and $a_1,\ldots,a_n$ be some real numbers
then, by Hoeffding's inequality,
%
%
\begin{equation}
\label{eqLaplace-Rademacher} \mathbb{E}\exp\Biggl(\sum_{i=1}^n
\varepsilon_i a_i \Biggr)\leq\exp\Biggl(
\frac{1}{2}\sum_{i=1}^n
a_i^2 \Biggr).
\end{equation}
We will also use a slightly modified version of the symmetrization
inequality: let $\mathcal{F}$ be a function class, $A_f,f\in\mathcal
{F}$ be a given function on $\mathcal{F}$ and $\Phi$ be a convex
nondecreasing function then
%
%
\begin{equation}
\label{eqsymmetrization-version2} \mathbb{E}\Phi\Bigl(\sup_{f\in
\mathcal{F}}
[Pf-P_n f-A_f ] \Bigr)\leq\mathbb{E}\Phi\Bigl(2\sup
_{f\in\mathcal{F}} [P_{n,\varepsilon
}f-A_f ] \Bigr),
\end{equation}
where $P$ is a measure, $P_n$ its associated empirical measure and
$P_{n,\varepsilon}$ the symmetrized empirical measure defined by
\[
P f=\mathbb{E}f(Z),\qquad P_nf=\frac{1}{n}\sum
_{i=1}^n f(Z_i)\quad\mbox{and}\quad
P_{n,\varepsilon}f=\frac{1}{n}\sum_{i=1}^n
\varepsilon_i f(Z_i),
\]
where $Z,Z_1,\ldots,Z_n$ are i.i.d. random variables distributed
according to $P$ and $\varepsilon_1,\ldots,\varepsilon_n$ are
independent Rademacher independent of $Z,Z_1,\ldots,Z_n$. The proof of
(\ref{eqsymmetrization-version2}) follows the same line as the
symmetrization inequality (cf. e.g., Theorem~2.1 in \cite{Kol11}).

Our analysis also relies upon some geometric argument. Indeed, the
strong convexity of the loss function in Assumption~\ref{assloss}
implies the $2$-convexity of the risk in the sense of \cite
{BarJorMcA06}; cf. (\ref{eqPisier}) for an explicit definition of the
$2$-convexity of a function $R(\cdot)$. This translates into a lower
bound on the gain obtained when applying Jensen's inequality to the
risk function $R$.
%
%
\begin{prop}\label{propstrong-2convexe}
Let $(X,Y)$ be a random couple in $\mathcal{X}\times\mathbb{R}$ and
$\mathcal{F}=\{f_1,\ldots,f_M\}$ be a finite dictionary in $
L_2(\mathcal{X},P_X)$ such that $|f_j(X)|\leq b,\ \forall j=1,\ldots,M$
and $|Y|\leq b$ a.s. Assume that, almost surely, the function $\ell
(Y, \cdot)$ is strongly convex with modulus of strong convexity
$C_\ell$ on $[-b,b]$. Then it holds that, for any $\theta\in\Theta$,
%
%
\begin{equation}
\label{eqPisier} R \Biggl(\sum_{j=1}^M
\theta_j f_j \Biggr)\leq\sum
_{j=1}^M\theta_j R(f_j)-
\frac{C_\ell}{2}\sum_{j=1}^M
\theta_j \Biggl\|f_j-\sum_{j=1}^M
\theta_j f_j \Biggr\|^2_{2}.
\end{equation}
\end{prop}
\begin{pf}
Define the random function $\ell(\cdot)=\ell(Y, \cdot)$. By strong
convexity and~\cite{HirLem01}, Theorem 6.1.2, it holds almost surely
that for any $a, a'$ in $[-b,b]$,
\[
\ell(a)\ge\ell\bigl(a'\bigr) + \bigl(a-a'\bigr)
\ell'\bigl(a'\bigr) + \frac{C_\ell}{2}
\bigl(a-a'\bigr)^2
\]
for any $\ell'(a')$ in the subdifferential of $\ell$ at $a'$.
Plugging $a=f_j(X)$, $a'=f_\theta(X)$, we get almost surely
\begin{eqnarray*}
&& \ell\bigl(Y, f_j(X)\bigr)
\\
&&\qquad \ge\ell\bigl(Y, f_\theta(X)\bigr)
+ \bigl(f_j(X)-f_\theta(X)\bigr)\ell'
\bigl(f_\theta(X)\bigr) +\frac{C_\ell}{2} \bigl[f_j(X)-f_\theta(X)
\bigr]^2.
\end{eqnarray*}
Now, multiplying both sides by $\theta_j$ and summing over $j$, we get
almost surely,
\[
\sum_j\theta_j\ell\bigl(Y,
f_j(X)\bigr)\ge\ell\bigl(Y, f_\theta(X)\bigr) +
\frac
{C_\ell}{2}\sum_j\theta_j
\bigl[f_j(X)-f_\theta(X)\bigr]^2.
\]
To complete the proof, it remains to take the expectation.
\end{pf}

\section{\texorpdfstring{Proof of Theorem~\protect\ref{THMAINUB}}
{Proof of Theorem A}}\label{secproofTheoA}

Let $x>0$ and assume that Assumption~\ref{assloss} holds throughout
this section. We start with some notation. For any $\theta\in\Theta
$, define
\[
\ell_\theta(y,x)=\ell\bigl(y,f_\theta(x)\bigr)\quad\mbox{and}\quad
\mathsf{R}(\theta)=\mathbb{E} \ell_\theta(Y,X)=\mathbb{E}\ell
\bigl(Y,f_\theta(X)\bigr),
\]
where we recall that $f_\theta=\sum_{j=1}^M\theta_j f_j$ for any
$\theta\in\mathbb{R}^M$. Let $0< \nu<1$. Let $(e_1,\ldots,e_M)$ is the
canonical basis of $\mathbb{R}^M$ and for any $\theta\in\mathbb
{R}^M$ define
\[
\tilde\ell_\theta(y,x)=(1-\nu)\ell_\theta(y,x)+\nu\sum
_{j=1}^M\theta_j\ell_{e_j}(y,x)\quad\mbox{and}\quad \tilde\mathsf{R}(\theta)=\mathbb{E} \tilde\ell_\theta(Y,X).
\]
We also consider the functions
\[
\theta\in\mathbb{R}^M\mapsto K(\theta)=\sum
_{j=1}^M\theta_j\log\biggl(
\frac{1}{\pi_j} \biggr)
\]
and
\[
\theta\in\mathbb{R}^M\mapsto V(\theta)=\sum
_{j=1}^M\theta_j\llVert
f_j-f_\theta\rrVert_2^2.
\]
Let $\mu>0$. Consider any oracle $\theta^*\in\Theta$ such that
\[
\theta^*\in\mathop{\operatorname{argmin}}_{\theta\in\Theta} \biggl(
\tilde
\mathsf{R}(\theta)+\mu V(\theta)+\frac{\beta}{n}K(\theta) \biggr).
\]

We start with a geometrical aspect of the problem. The two following
results follow from the strong convexity of the loss function $\ell$.

%
\begin{prop}\label{propconvexity}
When $\mu\leq(1-\nu)C_\ell/2$, the function $\theta\mapsto
H(\theta)=\tilde\mathsf{R}(\theta)+\mu V(\theta)+(\beta
/n)K(\theta)$ is
convex over the convex set $\Theta$.
\end{prop}
\begin{pf}
Let $\theta,\beta\in\Theta$ and $0\leq\alpha\leq1$. It follows
from some computation that
\[
V\bigl(\alpha\theta+(1-\alpha)\beta\bigr)=(1-\alpha)V(\beta)+\alpha
V(\theta)+
\alpha(1-\alpha)\llVert f_\theta-f_\beta\rrVert
_2^2.
\]
It follows from the strong convexity of $\ell(y,\cdot)$ that
\[
R\bigl(\alpha\theta+(1-\alpha)\beta\bigr)\leq(1-\alpha)R(\beta)+\alpha R(
\theta)-\frac{C_\ell}{2}\alpha(1-\alpha)\llVert f_\theta
-f_\beta\rrVert_2^2.
\]
Therefore, when $\mu\leq(1-\nu)C_\ell/2$, we have
\[
H\bigl(\alpha\theta+(1-\alpha)\beta\bigr)\leq(1-\alpha)H(\beta)+\alpha
H(\theta).
\]\upqed
\end{pf}

%
%
\begin{prop}\label{proptaylor}
Let $\mu\leq(1-\nu)C_\ell/2$. For any $\theta\in\Theta$,
\begin{eqnarray*}
\hspace*{-3pt}&& \tilde\mathsf{R}(\theta)-\tilde\mathsf{R}\bigl(\theta^*\bigr)
\\
\hspace*{-3pt}&&\qquad \geq\mu
\bigl(V
\bigl(\theta^*\bigr)-V(\theta) \bigr)+\frac{\beta}{n} \bigl(K\bigl
(\theta^*
\bigr)-K(\theta) \bigr)+ \biggl(\frac{ (1-\nu)C_\ell}{2}-\mu\biggr
)\llVert
f_\theta-f_{\theta^*}\rrVert_2^2.
\end{eqnarray*}
\end{prop}
\begin{pf}
Since $\theta\mapsto H(\theta)=\tilde\mathsf{R}(\theta)+\mu
V(\theta
)+(\beta/n)K(\theta)$ is convex over the convex set $\Theta$ and
$\theta^*$ is a minimizer of $H$ over $\Theta$, then there exists a
subgradient $\nabla H(\theta^*)$ such that for any $\theta\in\Theta
$ it holds, $\langle\nabla H(\theta^*),\theta-\theta^* \rangle\geq
0$. It yields
%
%
\begin{eqnarray}
\label{eqTaylor1}
&& \bigl\langle\nabla\tilde\mathsf{R}\bigl(\theta^*
\bigr), \theta-\theta^* \bigr\rangle\nonumber
\\
&&\qquad \geq\mu\bigl\langle\nabla V\bigl(\theta
^*\bigr), \theta^*-\theta \bigr\rangle+(\beta/n)\bigl\langle\nabla K\bigl(\theta^*\bigr
),\theta^*- \theta\bigr\rangle
\\
&&\qquad = \mu\bigl(V\bigl(\theta^*\bigr)-V(\theta) \bigr)-\mu\llVert f_\theta
-f_{\theta^*}\rrVert_2^2+(\beta/n) \bigl(K\bigl(
\theta^*\bigr)-K(\theta) \bigr).\nonumber
\end{eqnarray}
It follows from the strong convexity of $\ell(y, \cdot)$ that
\begin{eqnarray*}
\hspace*{-3pt}&& \tilde\mathsf{R}(\theta)-\tilde\mathsf{R}\bigl(\theta^*\bigr)
\\
\hspace*{-3pt}&&\qquad \geq\bigl\langle \nabla\tilde\mathsf{R}\bigl(\theta^*\bigr),\theta-\theta^* \bigr\rangle+
\frac{(1-\nu) C_\ell}{2}\| f_\theta-f_{\theta^*}\|_2^2
\\
\hspace*{-3pt} &&\qquad \geq\mu\bigl(V\bigl(\theta^*\bigr)-V(\theta) \bigr)+
\frac{\beta}{n} \bigl(K\bigl(\theta^*\bigr)-K(\theta) \bigr)+ \biggl(
\frac{(1-\nu) C_\ell}{2}-\mu\biggr)\llVert f_\theta-f_{\theta^*}\rrVert
_2^2,
\end{eqnarray*}
where the second inequality follows from the previous display.
\end{pf}

Let $\mathbf{H}$ be the $M\times M$ matrix with entries $\mathbf
{H}_{j,k}=\llVert f_j-f_{k}\rrVert_2^2$ for all $1\leq j,k\leq M$.
Let $s$ and $x$ be positive numbers and consider the random variable
\[
Z_n=(P-P_n) (\tilde\ell_{\hat\theta}-\tilde
\ell_{\theta^*})- \mu\sum_{j=1}^M\hat
\theta_j\llVert f_j-f_{\theta^*}\rrVert
^2_2-\mu\hat\theta\mathbf{H}\theta^*-
\frac{1}{s}K(\hat\theta).
\]
%

%
\begin{prop}\label{proprisk-bound}
Assume that
$10\mu\le\min(1-\nu, \nu)C_\ell$ and
$\beta\ge2n/s$. Then it holds
\[
R(\hat\theta)\leq\min_{1\leq j\leq M} \biggl[R(e_j)+
\frac{\beta
}{n}\log\biggl(\frac{1}{\pi_j} \biggr) \biggr]+ 2Z_n.
\]
\end{prop}
\begin{pf}
First note that the following equalities hold:
%
%
\begin{equation}
\label{eqeq1} \sum_{j=1}^M\hat
\theta_j\llVert f_j-f_{\theta^*}\rrVert
_2^2=V(\hat\theta)+\llVert f_{\hat\theta}-f_{\theta^*}
\rrVert_2^2
\end{equation}
and
%
%
\begin{equation}
\label{eqeq2} \hat\theta\mathbf{H}\theta^*=V(\hat\theta)+V\bigl(\theta
^*\bigr)+
\llVert f_{\theta^*}-f_{\hat\theta}\rrVert_2^2.
\end{equation}

It follows from the definition of $\hat\theta$ that
%
%
\begin{equation}
\label{eqeq-def-estimator} \tilde\mathsf{R}(\hat\theta)-\tilde\mathsf
{R}\bigl(\theta^*\bigr)
\leq(P-P_n) (\tilde\ell_{\hat\theta}-\tilde\ell_{\theta^*})+
\frac{\beta}{n} \bigl(K\bigl(\theta^*\bigr)-K(\hat\theta) \bigr).
\end{equation}
Then we use (\ref{eqeq1}) and (\ref{eqeq2}) together with (\ref
{eqeq-def-estimator}) to get
\begin{eqnarray}
\label{eqonevent1}
\nonumber
\tilde\mathsf{R}(\hat\theta)-\tilde\mathsf{R}\bigl(
\theta^*\bigr) 
&\le& 2\mu V(\hat\theta)+\mu V
\bigl(\theta^*\bigr)+2\mu\llVert f_{\hat\theta
}-f_{\theta^*}\rrVert
_2^2
\nonumber\\[-8pt]\\[-8pt]
&&{}+\frac{1}{s}K(\hat\theta)+\frac{\beta}{n} \bigl(K\bigl(\theta^*
\bigr)-K(\hat\theta) \bigr)+Z_n.\nonumber
\end{eqnarray}
Together with Proposition~\ref{proptaylor}, it yields
\[
\biggl(\frac{(1-\nu)C_\ell}{2}-3\mu\biggr)\llVert f_{\hat\theta
}-f_{\theta^*}
\rrVert_2^2\leq3\mu V(\hat\theta)+\frac
{1}{s}K(
\hat\theta)+Z_n.
\]
We plug the above inequality into (\ref{eqonevent1}) to obtain
\begin{eqnarray}
\label{eqtheo-final1}
\nonumber
\tilde\mathsf{R}(\hat\theta)-\tilde\mathsf{R}\bigl(\theta
^*\bigr)&\leq&\biggl(1+\frac{2\mu}{(1-\nu)C_\ell/2-3\mu} \biggr) \biggl
(\frac{1}{s}K(\hat
\theta)+Z_n \biggr)
\\
&&{}+\frac{\beta}{n} \bigl(K\bigl(\theta^*\bigr)-K(\hat\theta) \bigr)+\mu
V\bigl(\theta^*\bigr)
\nonumber
\\
&&{}+ \biggl(2\mu+\frac{6\mu^2}{(1-\nu)C_\ell/2-3\mu} \biggr)V(\hat\theta).\nonumber
\end{eqnarray}
Thanks to the $2$-convexity of the risk (cf. Proposition~\ref{propstrong-2convexe}), we have
\[
\tilde\mathsf{R}(\hat\theta)\geq\mathsf{R}(\hat\theta)+\nu(C_\ell/2)
V(\hat\theta).
\]
Therefore, it follows from (\ref{eqtheo-final1}) that
%
%
\begin{eqnarray}
\nonumber
\mathsf{R}(\hat\theta)&\leq&\tilde\mathsf{R}\bigl(\theta^*\bigr)+\mu V
\bigl(\theta^*\bigr)+\frac{\beta}{n}K\bigl(\theta^*\bigr)+ \biggl(1+
\frac{4\mu}{(1-\nu)C_\ell
-6\mu} \biggr)Z_n
\\
&&{} + \biggl(2\mu+\frac{12\mu^2}{(1-\nu)C_\ell-6\mu}-\nu\frac
{C_\ell}{2} \biggr) V(\hat\theta)
\\
&&{} + \biggl(\frac{1}{s}+\frac{4\mu}{s((1-\nu)C_\ell-6\mu
)}-\frac{\beta}{n} \biggr)K(
\hat\theta).\nonumber
\end{eqnarray}
Note now that $10 \mu\le\min(\nu, 1-\nu) C_\ell$ implies that
\[
\frac{4\mu}{(1-\nu)C_\ell-6\mu} \le1
\]
and
\[
2\mu+\frac
{12\mu^2}{(1-\nu)C_\ell-6\mu}-\nu
\frac{C_\ell}{2}\le0.
\]
Moreover, together, the two conditions of the proposition yield
\[
\frac{1}{s}+\frac{4\mu}{s((1-\nu)C_\ell-6\mu)}-\frac{\beta
}{n}\le0.
\]
%
Therefore, it follows from the above three displays that
\begin{eqnarray*}
\mathsf{R}(\hat\theta)&\leq&\min_{\theta\in\Theta} \biggl[\tilde\mathsf{R}
(\theta)+\mu V(\theta)+\frac{\beta}{n}K(\theta) \biggr]+2Z_n
\\
&\leq&\min_{j=1,\ldots,M} \biggl[\mathsf{R}(e_j)+
\frac{\beta}{n}\log\biggl(\frac{1}{\pi_j} \biggr) \biggr]+ 2Z_n.
\end{eqnarray*}\upqed
\end{pf}

To complete our proof, it remains to prove that $\mathbf{P}[Z_n >
(\beta
x)/n]\leq\exp(-x)$ and $\mathbb{E}[Z_n]\le0$ under suitable
conditions on
$\mu$ and $\beta$. Using, respectively, a~Chernoff bound and Jensen's
inequality, it is easy to see that both conditions follow if we prove
that $\mathbb{E}\exp( nZ_n/\beta) \le1$.
It follows from the excess loss decomposition:
\[
\tilde\ell_{\hat\theta}(y,x)-\tilde\ell_{\theta^*}(y,x)=(1-\nu) \bigl(
\ell_{\hat\theta}(y,x)-\ell_{\theta^*}(y,x)\bigr)+\nu\sum
_{j=1}^M\bigl(\hat\theta_j-
\theta_j^*\bigr)\ell_{e_j}(y,x)
\]
and the Cauchy--Schwarz inequality implies that it is enough to prove that
%
%
\begin{eqnarray}\label{EQconcentration-part11}
&& \mathbb{E}\exp\Biggl[s \Biggl((1-\nu)
(P-P_n) (\ell_{\hat\theta}-\ell_{\theta^*})
\nonumber\\[-8pt]\\[-8pt]
&&\hspace*{40pt}{} - \mu\sum_{j=1}^M
\hat\theta_j\llVert f_j-f_{\theta^*}\rrVert
^2_2-\frac{1}{s}K(\hat\theta) \Biggr) \Biggr]
\leq1\nonumber
\end{eqnarray}
and
%
%
\begin{equation}
\label{EQconcentration-part12}\qquad  \mathbb{E}\exp\Biggl[s \Biggl(\nu
(P-P_n) \Biggl(\sum
_{j=1}^M\bigl(\hat\theta_j-
\theta_j^*\bigr)\ell_{e_j} \Biggr)- \mu\hat\theta\mathbf{H}
\theta^*-\frac{1}{s}K(\hat\theta) \Biggr) \Biggr]\leq1
\end{equation}
for some $s \ge2n/\beta$ and assume this condition holds in the rest
of the proof.

%
We begin by proving~(\ref{EQconcentration-part11}). To that end,
define the symmetrized empirical process by $h\mapsto P_{n,\varepsilon
}h=n^{-1}\sum_{i=1}^n\varepsilon_ih(Y_i,X_i)$ where $\varepsilon
_1,\ldots,\varepsilon_n$ are $n$ i.i.d. Rademacher random variables
independent of the $(X_i,Y_i)$'s. Moreover, take $s$~and~$\mu$ such that
%
%
\begin{equation}
\label{eqcondition0} s\le\frac{\mu n}{[2C_b(1-\nu)]^2}.
\end{equation}
It yields
%
%
\begin{eqnarray}
\qquad &&\mathbb{E}\exp\Biggl[s \Biggl((1-\nu) (P-P_n) (
\ell_{\hat\theta}-\ell_{\theta^*})- \mu\sum_{j=1}^M
\hat\theta_j\llVert f_j-f_{\theta
^*}\rrVert
^2_2-\frac{1}s K(\hat\theta) \Biggr) \Biggr]
\nonumber
\\
&&\qquad \leq\mathbb{E}\exp\Biggl[s\max_{\theta\in\Theta} \Biggl((1-\nu)
(P-P_n) (\ell_{ \theta}-\ell_{\theta^*})\nonumber
\\
&&\hspace*{94pt}{}-\mu\sum_{j=1}^M\theta_j\llVert
f_j-f_{\theta^*}\rrVert^2_2-
\frac{1}sK(\theta) \Biggr) \Biggr]
\nonumber
\\
&&\qquad\leq\mathbb{E}\exp\Biggl[s\max_{\theta\in\Theta} \Biggl(2(1-\nu
)P_{n,\varepsilon}(\ell_{ \theta}-\ell_{\theta^*})
\nonumber\\[-8pt]\label{alignsym} \\[-8pt]
&&\hspace*{93pt}{} - \mu\sum
_{j=1}^M\theta_j\llVert
f_j-f_{\theta^*}\rrVert^2_2-
\frac{1}sK(\theta) \Biggr) \Biggr] \nonumber
\\
&&\qquad\leq\mathbb{E}\exp\Biggl[s\max_{\theta\in\Theta} \Biggl(2C_b
(1-\nu)P_{n,\varepsilon}(f_{\theta}-f_{\theta^*})
\nonumber\\[-8pt]\label{aligncontrac} \\[-8pt]
&&\hspace*{93pt}{} - \mu\sum_{j=1}^M\theta_j\llVert
f_j-f_{\theta^*}\rrVert^2_2-
\frac{1}sK(\theta) \Biggr) \Biggr],\nonumber
\end{eqnarray}
where~(\ref{alignsym}) follows from the slightly modified version of
the symmetrization inequality in (\ref{eqsymmetrization-version2})
and (\ref{aligncontrac}) follows from the contraction principle~\cite{LedTal91}, Theorem~4.12, applied to contractions
\[
\varphi_i(t_i)=C_b^{-1}\bigl[\ell
\bigl(Y_i, f_{\theta^*}(X_i)-t_i\bigr)-
\ell\bigl(Y_i, f_{\theta^*}(X_i)\bigr)\bigr]
\]
and $T\subset\mathbb{R}^n$ is defined by
\[
T=\bigl\{t \in\mathbb{R}^n\dvtx  t_i=f_{\theta^*}(X_i)-f_\theta(X_i),
\theta\in\Theta\bigr\}.
\]
Next, using the fact that the maximum of a linear function over a
polytope is attained at a vertex, we get
%
%
\begin{eqnarray}
\qquad &&\mathbb{E}\exp\Biggl[s \Biggl((1-\nu) (P-P_n) (
\ell_{\hat\theta}-\ell_{\theta^*})- \mu\sum_{j=1}^M
\hat\theta_j\llVert f_j-f_{\theta
^*}\rrVert
^2_2-\frac{1}s K(\hat\theta) \Biggr) \Biggr]
\nonumber
\\
&&\qquad\leq\sum_{k=1}^M\pi_k
\mathbb{E}\mathbb{E}_\varepsilon\exp\bigl[s \bigl(2C_b (1-
\nu)P_{n,\varepsilon}(f_k-f_{\theta^*})- \mu\llVert
f_k-f_{\theta^*}\rrVert^2_2 \bigr)
\bigr]
\nonumber
\\
&&\qquad\leq\sum_{k=1}^M
\pi_k\mathbb{E}\exp\biggl[\frac{[2C_b (1-\nu
)s)]^2}{2n}
\nonumber\\[-8pt] \label{alignrademacher}\\[-8pt]
&&\hspace*{89pt}{}\times \biggl(P_n-
\frac{2\mu n}{[2C_b(1-\nu)]^2s}P \biggr) (f_k-f_{\theta^*})^2
\biggr]\nonumber
\\
&&\qquad\leq\sum_{k=1}^M\pi_k
\mathbb{E}\exp\biggl[\frac{(2C_b(1-\nu
)s)^2}{2n}
\nonumber\\[-8pt]\label{aligncond0} \\[-8pt]
&&\hspace*{89pt}{}\times\biggl((P_n-P)
(f_k-f_{\theta^*})^2 -\frac{1}{4b^2}P(f_k-f_{\theta^*})^4
\biggr) \biggr],\nonumber
\end{eqnarray}
where (\ref{alignrademacher}) follows from (\ref
{eqLaplace-Rademacher}) and~(\ref{aligncond0}) follows from~(\ref
{eqcondition0}).
Together with the above display, Proposition~\ref{propBernstein-shifted} yields~(\ref{EQconcentration-part11})
as long as
%
%
\begin{equation}
\label{eqcondition3} s<\frac{n}{2\sqrt{3}bC_b(1-\nu)}.
\end{equation}
We now prove~(\ref{EQconcentration-part12}). We have
\begin{eqnarray*}
&&\mathbb{E}\exp\Biggl[s \Biggl(\nu(P-P_n) \Biggl(\sum
_{j=1}^M\bigl(\hat\theta_j-
\theta_j^*\bigr)\ell_{e_j} \Biggr)- \mu\hat\theta\mathbf{H}
\theta^*-\frac{1}{s}K(\hat\theta) \Biggr) \Biggr]
\\
&&\qquad \leq\sum_{j=1}^M\theta^*_j
\sum_{k=1}^M \pi_k\mathbb{E}
\exp\bigl[s \bigl(\nu(P-P_n) (\ell_{e_k}-
\ell_{e_j})-\mu\llVert f_j-f_k\rrVert
_2^2 \bigr) \bigr]
\\
&&\qquad \leq\sum_{j=1}^M\theta^*_j
\sum_{k=1}^M \pi_k\mathbb{E}
\exp\biggl[s\nu\biggl((P-P_n) (\ell_{e_k}-
\ell_{e_j})-\frac{\mu}{\nu C_b^2}P(\ell_{e_j}-
\ell_{e_k})^2 \biggr) \biggr]\leq1,
\end{eqnarray*}
where the last inequality follows from Proposition~\ref{propBernstein-shifted} when
%
%
\begin{equation}
\label{eqcondition4} s <\frac{2\mu n}{ C_b\nu(\nu C_b+4\mu b)}.
\end{equation}
It is now straightforward to see that the conditions of
Proposition~\ref{proprisk-bound}, the ones of~(\ref{eqcondition0}),
(\ref{eqcondition3}) and (\ref{eqcondition4}) are fulfilled when
\[
s=\frac{2n}{\beta},\qquad\mu=\min(\nu, 1-\nu)\frac{C_\ell}{10}%
\]
and
\[
\beta> \max\biggl[ \frac{8C_b^2(1-\nu)^2}{\mu}, 4\sqrt{3}bC_b(1-\nu),
\frac{C_b \nu(\nu C_b+4\mu b)}{\mu} \biggr].
\]




\printaddresses

\end{document}